\newif\iflandscape
\newif\ifportrait
\typein[\lorp]%
{Typein "l" (for landscape, twocolumn) or "p" (for portrait, onecolumn)}
\if l\lorp \landscapetrue \else \portraittrue \fi
\newlength{\extralineskip}
\ifportrait
  \documentstyle[]{article}
\fi
\iflandscape
  \documentstyle[twocolumn]{article}
\fi

\def\bea{\begin{eqnarray}}
\def\eea{\end{eqnarray}}

\def\beq{\begin{equation}}
\def\eeq{\end{equation}}
\def\ba{\beq\new\begin{array}{c}}
\def\ea{\end{array}\eeq}
\def\be{\ba}
\def\ee{\ea}

\def\f{1\over}
\makeatletter
\newdimen\normalarrayskip              
\newdimen\minarrayskip                 
\normalarrayskip\baselineskip
\minarrayskip\jot
\newif\ifold             \oldtrue            \def\new{\oldfalse}
\def\arraymode{\ifold\relax\else\displaystyle\fi} 
\def\eqnumphantom{\phantom{(\theequation)}}     
\def\@arrayskip{\ifold\baselineskip\z@\lineskip\z@
     \else
     \baselineskip\minarrayskip\lineskip2\minarrayskip\fi}
\def\@arrayclassz{\ifcase \@lastchclass \@acolampacol \or
\@ampacol \or \or \or \@addamp \or
   \@acolampacol \or \@firstampfalse \@acol \fi
\edef\@preamble{\@preamble
  \ifcase \@chnum
     \hfil$\relax\arraymode\@sharp$\hfil
     \or $\relax\arraymode\@sharp$\hfil
     \or \hfil$\relax\arraymode\@sharp$\fi}}
\def\@array[#1]#2{\setbox\@arstrutbox=\hbox{\vrule
     height\arraystretch \ht\strutbox
     depth\arraystretch \dp\strutbox
     width\z@}\@mkpream{#2}\edef\@preamble{\halign \noexpand\@halignto
\bgroup \tabskip\z@ \@arstrut \@preamble \tabskip\z@ \cr}%
\let\@startpbox\@@startpbox \let\@endpbox\@@endpbox
  \if #1t\vtop \else \if#1b\vbox \else \vcenter \fi\fi
  \bgroup \let\par\relax
  \let\@sharp##\let\protect\relax
  \@arrayskip\@preamble}
\def\eqnarray{\stepcounter{equation}%
              \let\@currentlabel=\theequation
              \global\@eqnswtrue
              \global\@eqcnt\z@
              \tabskip\@centering
              \let\\=\@eqncr
              $$%
 \halign to \displaywidth\bgroup
    \eqnumphantom\@eqnsel\hskip\@centering
    $\displaystyle \tabskip\z@ {##}$%
    global\@eqcnt\@ne \hskip 2\arraycolsep
         $\displaystyle\arraymode{##}$\hfil
    global\@eqcnt\tw@ \hskip 2\arraycolsep
         $\displaystyle\tabskip\z@{##}$\hfil
         \tabskip\@centering {##}\tabskip\z@\cr}
\def\theequation{\thesection.\arabic{equation}}
\newfont{\hr}{msbm10}
\newfont{\ams}{msam10}
\begin{document}
\begin{titlepage}
\setcounter{footnote}0
\begin{center}
\hfill ITEP-M1/95\\
\hfill FIAN/TD-19/94\\
\hfill q-alg/9501013\\
\vspace{0.3in}
\ifportrait
{\LARGE\bf Non-Standard KP Evolution}\\{\LARGE\bf and Quantum $\tau$-function}
\fi
\iflandscape
{\LARGE\bf Non-Standard KP Evolution and Quantum $\tau$-function}
\fi
\\[.4in]
{\Large S. Kharchev\footnote{E-mail address:
kharchev@vxitep.itep.ru, kharchev@lpi.ac.ru}$^{\dag\ddag}$,
A.Mironov\footnote{E-mail address: mironov@lpi.ac.ru}$^{\ddag\dag}$,
A.Morozov\footnote{E-mail address:
morozov@vxdesy.desy.de}$^{\dag}$}\\
\bigskip $\phantom{gh}^{\dag}$  {\it ITEP, Moscow, 117 259, Russia}\bigskip\\
$\phantom{gh}^{\ddag}$  {\it Theory Department,  P. N. Lebedev
Physics Institute, Leninsky prospect, 53,\\ Moscow,~117924, Russia}
\end{center}
\bigskip
\bigskip

\begin{quotation}
One possible way to fix partly a ``canonical definition''
of $\tau$-functions beyond the conventional KP/Toda framework
could be to postulate that evolution operators are
{\it always} group elements. We discuss implications of this
postulate for the first non-trivial case: fundamental
representations of quantum groups $SL_q(N)$. It appears that
the most suited (simple) for quantum deformation framework
is some non-standard formulation of KP/Toda systems. It turns
out that the postulate needs to be slightly modified
to take into account that no ``nilpotent subgroups''
exist in $SL_q(N)$ for $q\neq 1$. This has some definite and
simple implications for $q$-determinant-like representations of
quantum  $\tau$-functions.
\end{quotation}
\end{titlepage}
\clearpage
\newpage
\section{Introduction}

Appropriately (and broadly) generalized classical integrable
hierarchies are now widely believed to describe non-perturbative
effective actions of string models while their quantum counterparts
should be relevant for description of the second-quantized string
theory. The modest purpose of these notes is to illustrate
some pecularities on the way from classical to quantum hierarchies
of the simplest -- KP/Toda -- type. One of the possible approaches
to quantization is to make use of the group theory
interpretation of hierarchies and then just substitute ordinary groups
by their quantum deformations. This is the line to be followed
in the present paper.

The basic object in the theory of hierarchies is $\tau$-function --
the generating functional of all the matrix elements of a given
group element $g \in G$ in a given (highest weight) representation
$R$:
\be
\tau_R (t,\bar t|g) \equiv \sum_{\{m,\bar m\}\in R}
s^{R}_{m,\bar m}(t,\bar t) <m |g| \bar m>
\label{taugen}
\ee
The choice of the functions $s^{R}_{m,\bar m}(t,\bar t)$
is the main ambiguity in the definition of $\tau$-function
and needs to be made in some clever way (not yet known
in full generality). In the case of the highest weight
representation $R$, it can be partly fixed by the requirement that
\be
\tau_R (t,\bar t|g)  =
<0_R| U(t) g \bar U(\bar t)| 0_R >
\label{tau}
\ee
where operators $U$ and $\bar U$ do not depend on $R$.

KP/Toda-type $\tau$-function arises when $G = SL(N)$ and $R$ is
one of the $N-1$ fundamental representations. All fundamental
representations of $SL(N)$ are skew products of the first
($N$-dimensional) one $F \equiv F_1$: $F_n = \wedge^n F$,
and thus can be also described in terms of (fermionic) intertwining
operators. Entire KP/Toda hierarchy is obtained in the limit of
$N \rightarrow \infty$ (with $n$ playing the role
of the ``zero-time''), and has also an alternative description
in terms of the level $k=1$ Kac-Moody algebras. We shall,
however, concentrate on the case of the generic $N$.

Operators $U(t)$, $\bar U(\bar t)$, when restricted onto any
fundamental representation $F_n$, turn into
\be
U(t) = \exp \left(\sum_{k\geq 1} t_kT_+^{(k)}\right), \\
\bar U(\bar t) = \exp \left(\sum_{k\geq 1} \bar t_kT_-^{(k)}\right),
\label{Uop}
\ee
where $\displaystyle{T_{\pm}^{(k)} = \sum_{\vec\alpha : h(\vec\alpha) = k}
T_{\pm\vec\alpha}}$ are sums of all the generators of $SL(N)$,
associated with the positive/negative roots of ``height'' $k$
(in the {\it first} fundamental representation $F$ such $T_{\pm}^{(k)}$
look like $N\times N$ matrices with units at the $k$-th upper/lower
diagonal and zeroes elsewhere). The crucial feature of these
operators is their commutativity:
\be
\left[ T_+^{(k)},\ T_+^{(l)} \right] = 0,\ \ \
\left[ T_-^{(k)},\ T_-^{(l)} \right] = 0.
\ee
It is a peculiarity of fundamental representations of $SL(N)$ that such
simple $U(t)$, $\bar U(\bar t)$ are sufficient to generate {\it all}
the elements of representation (generally $U$, $\bar U$ depend on more
time-variables and more mutually non-commuting generators of $G$).

Operators, defined in (\ref{Uop}) have the following properties:

(i) $U,\ \bar U\ \in\  G$;

(ii) more specific, $U$, $\bar U$ belong to the ``nilpotent subgroup''
$NG$ of $G$. Actualy $NG$ is a subgroup of Borel subgroup:
$NG \subset BG \subset G$ (in $F_1$ $BG$ consists of all the
upper triangular matrices with unit determinant, while
matrices from $NG$ are additionally constrained to have
only unit elements on the main diagonal);

(iii) if comultiplication is defined so that
\be
\Delta(T_{\pm\vec\alpha}) = T_{\pm\vec\alpha}\otimes I +
	 I \otimes T_{\pm\vec\alpha}
\label{claco}
\ee
then
\be
\Delta U(t) = U(t) \otimes U(t) =
\left(U(t)\otimes I\right) \left(I\otimes U(t)\right)
\label{clacoU}
\ee

Shortly speaking, evolution operators $U$, $\bar U$ are
just ``group elements'' of $NG$.\footnote{
  We refer to \cite{GKLMM} and \cite{MV,Mon} for a lengthy discussion
        of what we mean by ``group elelment'' in this context. In a
        word, this is what has been called ``universal $T$-matrix'' in the
quantum group theory \cite{FTR,FG,JS}.
}

These properties seem to be rather appealing and it is natural to try
to preserve them in any generalization. Below we consider a generalization
in one particular direction -- that of quantum groups: in what
follows we discuss $\tau$-functions for fundamental representations
of $SL_q(N)$. There are two immediate things to be taken into account.
First, there is nothing similar to the operators $T_{\pm}^{(k)}$
for $q\neq 1$ (at least, nothing what could be defined in terms of generators
without any references to a specific representation $R$). This implies that
explicit expressions for $U(t)$ and $\bar U(\bar t)$ should be very different
from (\ref{Uop}).  Second, there is no reasonable notion of the nilpotent
subgroup $NG_q$: only quantum deformation of the Borel subgroup $BG_q \subset
G_q$ is nicely defined.  Because of this one should not insist on validity of
eq.(\ref{clacoU}): it requires some modification.

According to \cite{MV} parametrization of the group elements, which admits
the most straightforward quantum deformation, involves only the
{\it simple} roots $\pm \vec\alpha_i$, $i = 1,\ldots,r_G$:
\be
g = g_Ug_Dg_L, \\
g_U =\left.\prod_s\right.^<\ e^{\theta_sT_{i(s)}},   \ \ \
g_L = \left.\prod_s\right.^>\ e^{\chi_sT_{-i(s)}}, \ \ \
g_D = \prod_{i=1}^{r_G} e^{\vec\phi\vec H}
\label{genpar}
\ee
Every particular simple root $\vec\alpha_i$ can appear several
times in the product, and there are different parametrizations of
group elements of such a type, depending on the choice of the
set $\{s\}$ and the mapping $i(s)$ of this set into that of simple roots.
Quantum deformation of {\it such} formula is especialy simple because
comultiplication rule is simple for the generators, associated with
the {\it simple} roots:  \be \Delta(T_i) = T_i\otimes q^{-2H_i} + I\otimes T_i,
\\ \Delta(T_{-i}) = T_{-i}\otimes I + q^{2H_i}\otimes T_{-i} \label{quco} \ee
For $q\neq 1$ any expression of the form (\ref{genpar})
remains just the same, provided exponentials in $g_U$ and $g_L$
are understood as $q$-exponentials (in the simply-laced case,
$q^{||\vec\alpha_i||^2/2}$-exponentials in general), and parameters
$\theta, \chi, \vec\phi$ become non-commuting generators of the
``coordinate ring'' of $G_q$. Actually, they form a kind of
Heisenberg algebra:
\be
\theta_s\theta_{s'} = q^{-\vec\alpha_{i(s)}\vec\alpha_{i(s')}}
\theta_{s'}\theta_s,\ \ \ s<s', \\
\chi_s\chi_{s'} = q^{-\vec\alpha_{i(s)}\vec\alpha_{i(s')}}
\chi_{s'}\chi_s,\ \ \ s<s', \\
e^{\vec\beta\vec\phi}\theta_s = q^{\vec\beta\vec\alpha_{i(s)}}
\theta_s e^{\vec\beta\vec\phi}, \\
e^{\vec\beta\vec\phi}\chi_s = q^{\vec\beta\vec\alpha_{i(s)}}
\chi_s e^{\vec\beta\vec\phi}
\label{comrel}
\ee
These relations imply that\footnote{See \cite{MV}
for all the notations and definitions.} $\Delta(g) = g\otimes g$.

The simplest possible assumption about evolution operators
would be to say that $U(t)$ is always an object of the type $g_U$,
while $\bar U(\bar t)$ -- of the type $g_L$. However, these are no
longer group elements:
$$
\Delta(g_U) \neq g_U \otimes g_U, \ \ \
\Delta(g_L) \neq g_L \otimes g_L,
$$
because of the lack of factors $g_D$. This is the exact meaning of the
claim that there is no ``nilpotent subgroup'' $NG_q$ (but $BG_q$
does exist, since $\Delta(g_Ug_D) = (g_Ug_D)\otimes(g_Ug_D)$).
Despite this ``problem'' we will insist on identification of
$U$ and $\bar U$ as objects of the type $g_U$ and $g_L$ respectively,
and will explicitly investigate implications of the failure of
(\ref{clacoU}) (see Conclusion where another, perhaps more attractive, option
is
mentioned).  In fact, instead of (\ref{clacoU}) we will have

\be
\Delta(U(\xi)) = U^{(2)}_L(\xi) \cdot
U^{(2)}_R(\xi),
\ee
where
\be\label{1.11}
U(\xi) = \left.\prod_s\right.^<\ {\cal
E}_q\left(\xi_sT_{i(s)}\right),
\ee
\be
U^{(2)}_L = \left.\prod_s\right.^<\ {\cal
E}_q\left(\xi_sT_{i(s)}\otimes q^{-2H_{i(s)}}\right) \neq I\otimes U(\xi), \\
U^{(2)}_R = \left.\prod_s\right.^<\
{\cal E}_q\left(\xi_s I\otimes T_{i(s)}\right) = I\otimes U(\xi)
\ee
and this will have simply accountable implications for determinant
formulas for quantum $\tau$-functions.

In what follows we first discuss various interesting ways to
specify the map $i(s)$ in the case of fundamental representations.
Among these there is especially simple one,
$s = 1,\ldots,r_G$, $i(s) = s$. However, it gives rise to $U(t)$
which is different from (\ref{Uop}) even in the classical case
of $q=1$. Therefore, we briefly describe the classical
hierarchy with this non-standard evolution.  Finally, we consider the
corresponding quantum
deformation and derive the substitute of the determinant formulas for $\tau_n
\equiv \tau_{F_n}$ in the case of $q \neq 1$.
Let us stress that by the
multiplication of the evolution operators $U$ and $\bar U$ and the group
element
$g$ in the definition of $\tau$-function (\ref{tau}) we always understand
the group
multiplication law, i.e. elements of algebra $\theta,\ \phi$ and $\chi$
(\ref{comrel}) in evolution operators commute with elements of
the corresponding algebra in $g$, see \cite{MV,Mon}. This is very essential
for the determinant formulas of section 4.

\section{Group elements through simple roots: examples}

We briefly discuss here three natural choices of parametrization
(\ref{genpar}) of the group elements.

As we already mentioned, every parametrization of this form is
straightforwardly deformed to $q\neq 1$ \cite{MV}.
The most economic way to parametrize in this way the entire
group manifold of $SL(N)$ is to take
$s = 1,\ldots,\frac{N(N-1)}{2}$ and the map
$$
i(s):\ \ \ 1,2,\ldots,r-1,r;\  1,2,\ldots,r-1;\  1,2,3;\  1,2;\  1;
\\ r = {\rm rank}\  SL(N) = N-1
$$
i.e.
\be
U(\xi) = \prod_{1\leq i \leq N}\prod_{i<j\le N}\exp\left(
\xi_{ij} T_{j-i}\right).
\label{i(s)}
\ee

This is, however, a little too much for our purposes. The orbits
of $SL(N)$ in {\it fundamental} representations can be
parametrized by only $r = N-1$ parameters, and the purpose is
to find an adequate parametrization of such submanifolds.
This is easy to do in the classical ($q=1$) case, and, at least,
three natural possibilities will be considered in this section.
However, of these three only one will be easily deformed, and
it is the one with no direct relation to conventional evolution
(\ref{Uop}).

{\bf Parametrization A.} The simplest possibility is just to
restrict the set $\{s\}$ to $s = 1,\ldots,r$ and take
$i(s) = s$, i.e. take
\be
U^{(A)}(\xi) = \left.\prod_{i=1}^{r_G}\right.^<\
\exp\left(\xi_iT_i\right).
\label{UA}
\ee
This is enough to generate all the states of any fundamental
representation from the corresponding vacuum (highest vector)
state, but
$<0_{F_n}|\ U^{(A)}(\xi) $ has little to do with
$<0_{F_n}|\ U(t)$  (where $U(t)$ is given by (\ref{Uop})).
It can be better to say that identification
$ <0_{F_n}|\ U^{(A)}(\xi) = <0_{F_n}|\ U(t)$
defines a relation $\xi_i(t)$, which explicitly depends on $n$.

One can build the theory of the KP/Toda hierarchies
in terms of $\xi$-variables instead of conventional $t$-variables
(see a brief discussion in s.3 below), but it can {\it not} be obtained
just by change of time-variables: the whole construction
looks different. For it, this new construction is immediately
deformed to the case of $q\neq 1$: instead of (\ref{UA}) we
just write
\be
U^{(A)}(\xi) = \left.\prod_{i=1}^{r_G}\right.^<\
{\cal E}_q\left(\xi_iT_i\right),
\label{UAq}
\ee
where $\xi$'s are non-commuting variables,
\be
\xi_i\xi_j = q^{-\vec\alpha_i\vec\alpha_j}\xi_j\xi_i,
\ \ \ i<j,
\ee
and it is easy to derive a quantum counterpart for any statement
of the classical ($q=1$) theory once it is formulated for the
$\xi$-parametrization (see s.4 for some results in the $q\neq 1$
case).

{\bf Parametrization B (conventional).} Of course, one can
insist on using the conventional $t$-variables, i.e. to make
the identification of the group elements
\be
U^{(B)}(\xi) = \prod_s \exp\left(\xi_sT_{i(s)}\right) =
U(t) = \exp \left(\sum_k t_kT_+^{(k)}\right)
\label{2.16}
\ee
(which implies that
$ <0_{F_n}|\ U^{(B)}(\xi) = <0_{F_n}|\ U(t)$
with some $n$-independent functions $\xi_s(t)$).
The difference between the two expressions in (\ref{2.16}) is that the r.h.s.
contains mutually commuting combinations of (non-simple) root generators, while
the l.h.s. contains only (mutually non-commuting) simple-root generators. Such
reparameterization indeed exists, but the set $\{s\}$ should contain, at least,
$\frac{N(N-1)}{2}$ elements and one can take $i(s)$ as in (\ref{i(s)}) -- the
only thing is that now not all of the $\xi_s$ are independent: instead they are
expressed through $r = N-1$ time-variables $t_k$.  For example, the
$t_1$-dependence of $\xi_{ij}$ is given by
\be
\xi_{ij} = \frac{t_1}{N+i-j} + {\cal O}(t_2, t_3,
\ldots).  \label{xivert}
\ee

{\it Open problem}.
In order to get a reasonable quantum deformation of parametrization B,
one needs to reproduce the proper commutation relations
\be
\xi_s\xi_{s'} = q^{-\vec\alpha_{i(s)}\vec\alpha_{i(s')}}
\xi_{s'}\xi_s,\ \ \ s<s',
\ee
between the $\frac{N(N-1)}{2}$ variables $\xi_s$
as a corollary of {\it some} relations between $r=N-1$
variables $t_k$ (which, of course, do not need to commute
when $q\neq 1$). To make this possible, one should also somehow
deform the relations (\ref{xivert}) at $q\neq 1$.
This is a separate problem, which we do not have
immediate solution to.

{\bf Parametrization C (Miwa variables).} One more option is
to remain with the conventional time-variables $t_k$, but make
the (representation-independent) Miwa transform
$t_k = \frac{1}{k}\sum_a \lambda_a^k$. This Miwa parametrization
is, in fact, perfectly consistent with the simple-root
decomposition:
\be
U(t) = \prod_a \exp\left(\sum_{k=1}^{r_G}\frac{\lambda_a^k}{k}
T_+^{(k)}\right) =
\prod_a \left(\prod_{i=1}^{r_G} e^{\lambda_aT_i}\right).
\label{claMipa}
\ee
The set $\{s\}$ and maping $i(s)$ here are not of the ``most
economic'' type (\ref{i(s)}), but the general rule (\ref{comrel})
of the quantum deformation is, of course, applicable.

{\it Open problem}. However, (\ref{i(s)}) implies
the quantum formula in the form
\be
\prod_a \left(\prod_{i=1}^{r_G} {\cal E}_q\left(\lambda_{ai}T_i
\right)\right),
\ee
where $\lambda_{ai}$ with different $i$ and the same $a$ do
not commute. At the same time, the constraint $\lambda_{ai} =
\lambda_{aj}$ for $i\ne j$ is of crucial importance for the classical ($q=1$)
formula (\ref{claMipa}). What is the proper deformation of this
constraint remains unclear. Solution of the puzzle should
probably exploit the fact that, in the classical case, the constraint
selects out irreducible representations of the coordinate ring
(dual algebra) of $G$.

\section{Classical ($q=1$) KP/Toda theory}
\subsection{Determinant formulas and systems of equations}
Let us begin with reminding the main features of the standard
construction of hierarchies, starting from the fundamental
representations of $SL(N)$ group. At the moment, we do not
specify any parameterization, working with some arbitrary evolution
$U(t)$.

We consider $\tau_n \equiv \tau_{F_n}(t,\bar t|g)$.
These $\tau$-functions
satisfy Hirota equations, which can be derived by the general
procedure from the action of intertwining operator (fermion)
$\psi^\pm :\ \ F_1\otimes F_n$~{\ams \symbol{29}}~$F_{n+1}$ (see
\cite{GKLMM} for details). Hirota equations, which are straightforwardly
deformed to the $q\neq 1$ case \cite{GKLMM}, are not
specific for the fundamental representations of $SL(N)$ and
for the KP/Toda $\tau$-functions. We will turn to these equations a little
later, but first discuss a more specific subject: determinant
representations of $\tau_n$ (these, at least in the simple form, are specific
for fundamental representations and/or level $k=1$ Kac-Moody algebras -- i.e.
the cases where {\it fermionization}, not only {\it bosonization}
is available).

Let us first consider
$$
\tau_1 = <0_{F_1}|U(t) g \bar U(\bar t)|0_{F_1}>.
$$
Note that a specific feature of $F=F_1$ is
$$
<0_{F}|U(t) =
\sum_k {\cal P}_k(t) <0_{F} | T_+^k =
\sum_k {\cal P}_k(t) <k_{F}|,
$$
where the r.h.s. is reexpanded in terms of generalized Schur polinomials
(the first equality in this formula defines these polynomials)
and the $N$ states of $F = F_1$ are denoted by
$<k_F| = <0_F| T_+^k$, $k = 0,...,r=N-1$.
Thus,
\be
\tau_1(t,\bar t|g) =
\sum_{k,\bar k} {\cal P}_k(t){\cal P}_{\bar k}(\bar t)
<0_{F}|T_+^kg T_-^{\bar k}|0_{F}> =\\=
\sum_{k,\bar k} {\cal P}_k(t){\cal P}_{\bar k}(\bar t)
<k_{F}| g |\bar k_{F}>=\sum_{k,\bar k} {\cal P}_k(t)g_{k,\bar k}{\cal P}_{\bar
k}(\bar t).  \label{clatau1}
\ee
One can also define
\be
\tau_1^{m\bar m}
\equiv <m_F|U(t) g \bar U(\bar t)|\bar m_F> = \sum_{k,\bar k} {\cal P}_k(t)
g_{m+k,\bar m+\bar k}{\cal P}_{\bar k}(\bar t).
\label{clataumm1}
\ee

Now we can turn to generic fundamental representation $F_n$.
Since
\be
\left.<m_1\ldots m_n\right._{F_n}| =
<{m_1}_F|\otimes<{m_2}_F|\otimes\ldots\otimes<{m_n}_F| +\\
+
{\rm antisymmetrization\ over\ }m_1,\ldots,m_n =\\
= \sum_P (-)^P <m_{P(1)}|\otimes <m_{P(2)}| \otimes\ldots\otimes
<m_{P(n-1)}|
\label{as}
\ee
the vacuum (highest weight) state of $F_n$ can be written as
\be
<0_{F_n}| = <0,1,\ldots,n-1_{F_n}|=\\
= \sum_P (-)^P <P(0)_F|\otimes <P(1)_F| \otimes\ldots\otimes
<P(n-1)_F|.
\label{ass}
\ee
Taking into account that
\be\label{coprodcl}
U(t)|_{F_n} = \Delta^{n-1}U(t) = U(t)^{\otimes n}, \ \ \
g|_{F_n} = \Delta^{n-1}(g) = g^{\otimes n}
\ee
one finally gets:
\be
\tau_{n+1}(t,\bar t|g) \equiv <0_{F_n}| U(t) g \bar U(\bar t)
|0_{F_n}> =  \\
= \sum_{P,\bar P} (-)^P(-)^{\bar P}
\prod_{k=0}^{n} <P(k)_F|U(t) g \bar U(\bar t)|\bar P(k)_F> = \\
= \det_{0\leq m,\bar m<n}\tau_1^{m\bar m} = \det_{0\leq m,\bar m<n}
\sum_{l,\bar l} {\cal P}_{l-m}(t)g_{l_F,\bar l_F}{\cal P}_{\bar l-\bar
m}(\bar t)=\\ = \sum_{{1<m_1<m_2<...}\atop{1<\bar m_1<\bar m_2<...}}
\det_{ji}{\cal P}_{m_j-i}(t) \det_{ji}g_{m_{j}\bar m_{i}} \det_{ij}
{\cal P}_{\bar m_i-\bar j}(\bar t).  \label{det}
\ee
These determinant
formulas are equivalent to the KP/Toda equations.
However, determinant formulas are not the
simplest starting point to derive the equations (see, for example, \cite{MMV}).
Here we apply the simpler method \cite{DJKM,GKLMM} making use of
fermionic intertwining operators $\psi_i^{\pm}$ ($i=1\ldots N$). The key
ingredient of the derivation (see \cite{GKLMM,Mon}) is the composite
intertwining operator $\Gamma=\sum_i\psi_i^+\otimes\psi_i^-$.
In order to
get a set of equations, one consider the matrix element of the interwiner
identity $\Gamma(g\otimes g)=(g\otimes g)\Gamma$ between the states
$\langle 0_{F_{n+1}}|U(t)\otimes \langle 0_{F_{m-1}}|U(t')$ and
$\bar U(\bar t)|0_{F_{n}}\rangle\otimes \bar U(\bar t')|0_{F_{m}}\rangle$.

\be\label{BIglinfty}
\sum_i \langle 0_{F_{n+1}} | U(t) \psi_i^+ g \bar U(\bar t)
|0_{F_{n}}\rangle \cdot \langle 0_{F_{m-1}} | U(t')
\psi_i^{-} g \bar U(\bar t')
|0_{F_{m}}\rangle =\\
=\sum_i \langle 0_{F_{n+1}} | U(t) g\psi_i^+  \bar U(\bar t)
|0_{F_{n}}\rangle \cdot \langle 0_{F_{m-1}} | U(t')
g\psi_i^{-}  \bar U(\bar t')
|0_{F_{m}}\rangle.
\ee
One can rewrite (\ref{BIglinfty}) through the free fermion fields
$\psi^+ (z)\equiv \sum_{i=1}^N \psi_i^+ z^i$ and $\psi^{-}(z)\equiv\sum_{i=1}^N
\psi_i^{-}z^{N-i+1}$:
\be
\oint_{\infty} {dz\over z^{N+2}}
\langle 0_{F_{n+1}} | U(t) \psi^+(z) g \bar U(\bar t)
|0_{F_{n}}\rangle \cdot \langle 0_{F_{m-1}} | U(t')
\psi^{-}(z) g \bar U(\bar t')
|0_{F_{m}}\rangle =\\=
\oint_0{dz\over z^{N+2}}
\langle 0_{F_{n+1}} | U(t) g\psi^+(z)  \bar U(\bar t)
|0_{F_{n}}\rangle \cdot \langle 0_{F_{m-1}} | U(t')
g\psi^{-}(z)  \bar U(\bar t')
|0_{F_{m}}\rangle.
\ee
The same formulas can be written in more compact form using the
following (Baker-Akhiezer) functions

\be\label{BA}
\Psi^{\pm,i}_n\equiv \langle 0_{F_{n\pm 1}} | \hat U(t)
\psi_i^{\pm} g \hat{\bar U}(\bar t).
|0_{F_{n}}\rangle
\ee
Then (\ref{BIglinfty}) can be rewritten
in the form
\be\label{BIA}
\sum_i\Psi^{+,i}_k(t)\Psi^{-,i}_l(t')=
\sum_j\bar\Psi^{+,j}_{k+1}(\bar t)\bar\Psi^{+,j}_{l-1}(\bar t'),
\ee
where $\bar\Psi$ is defined analogously to (\ref{BA}) but with the fermion
situated to the right of the group element $g$.

One can also define vertex operators which generates Baker-Akhiezer functions:

\be
\sum_i\Psi^{+,i}_k(t)z^{i}\equiv \hat X^{+} (z,t)\tau_n(t),\ \ \
\sum_i\Psi^{-,i}_k(t)z^{N-i+1}\equiv \hat X^{-} (z,t)\tau_n(t),
\label{VO}
\ee
and analogously for ${\hat {\bar X}}^{\pm}(z,t)$.
Then, (\ref{BIglinfty}) can be also written as

\be\label{IntHir}
\oint{dz\over z^{N+2}} \hat X^-(z,t)\tau_n(t,\bar t) \hat
X^{+}(z,t')\tau_m(t',\bar t')=\\= \oint{dz\over z^{N+2}}\hat
{\bar X}^-(z,\bar t)\tau_{n+1}(t,\bar t) \hat {\bar X}^{+}(z,\bar
t')\tau_{m-1}(t',\bar t'),
\ee

Now we apply all these formulas to the particular choices of evolution
operators $U(t)$, $\bar U(\bar t)$.

\subsection{Conventional parametrization (B)}

This evolution leads to the standard KP/Toda hierarchy.
One makes use of expressions (\ref{Uop}). It gives
$$
<0_{F}|U(t) = <0_{F} | \exp\left(\sum_k t_kT_+^{(k)}\right)=
\sum_k P_k(t) <k_{F}|,
$$
with the orthodox Schur polynomials $P_k(t)$ defined by
$\exp\left(\sum_k t_kz^k\right) = \sum_k P_k(t)z^k$.
The main peculiarity of this evolution is the property
\be
\tau_1^{m\bar m} = \frac{\partial}{\partial t_m}
\frac{\partial}{\partial\bar t_{\bar m}}\tau_1 =
\left(\frac{\partial}{\partial t_1}\right)^m
\left(\frac{\partial}{\partial\bar t_1}\right)^{\bar m}\tau_1.
\label{clat-derB}
\ee

Determinant formula (\ref{det}) in this parametrization
\be
\tau_n(t,\bar t|g)=
\sum_{{1<m_1<m_2<...}\atop{1<\bar m_1<\bar m_2<...}} \det_{ji}P_{m_j-i}(t)
\det_{ji}g_{m_{j}\bar m_{i}}
\det_{ij} P_{\bar m_i-\bar j}(\bar t).  \label{detB}
\ee
leads to equations which are nothing but
bilinear Pl\"ucker relations \cite{Sat}, while (\ref{detB})
demonstrates explicitly that $\tau$-function is spanned by the particular
Pl\"ucker coordinates -- Schur functions $\xi_=\det_{ij}P_{m_j-i}(t)$, where
$$ labels oung tables -- see \cite{Mac}).

In order to obtain the equations in parametrization B
one should note that, in this case, vertex operators (\ref{VO}) are
\be
\hat X^{+}(z,t) = \hbox{Pr}_N\left[e^{\xi(z,t)}
\hbox{Pr}\left[z^ne^{-\xi(z^{-1},\tilde
\partial_{t})}\right]\right],\\
\hat X^{-}(z,t) = \hbox{Pr}_N\left[e^{-\xi(z,t)}\hbox{Pr}\left[z^{N-n+1}
e^{\xi(z^{-1},\tilde
\partial_{t})}\right]\right],\\
\ \ \ \xi(z,t)\equiv\sum_i^Nz^it_i,\ \ \tilde
\partial_{t_k}\equiv {\f k} \partial_{t_k}
\ee
(and similarly for the right
vacuum state), where $\hbox{Pr}[f(z)]$ projects onto the polynomial part of
function
$f(z)$ and $\hbox{Pr}_l[f(z)]$ projects onto the polynomial part of the
degree $l$.

One can explicitly write equations (\ref{IntHir}) in the integral form,
which can be easily transformed to an infinite
set of differential equations by expanding in powers of time differences
$t_i-t_i'$ etc. \cite{DJKM}. As $N\longrightarrow\infty$, these equations
give rise to the standard Toda Lattice hierarchy.

In particular, the simplest equation, which contains only derivatives with
respect to the first times, is

\be
{\partial\tau_n(t,\bar t)\over\partial t_1}\cdot
{\partial\tau_n(t,\bar t)\over\partial \bar t_1}-\tau_n(t,\bar t)\cdot
{\partial^2\tau_n(t,\bar t)\over\partial t_1 \partial\bar t_1}=
\tau_{n+1}(t,\bar t)\tau_{n-1}(t,\bar t).
\label{1Hireq}
\ee
Equation of such a type can be also obtained relatively easy
from the determinant formulas \cite{MMV}.

\subsection{On KP/Toda hierarchy in parametrization A}

Now  we consider the same conventional hierarchy with a different evolution A.
Our purpose is to demonstrate that in this
parametrization the main features of the hierarchy
are preserved -- there are determinant formulas and a hierarchy of differential
equations.

{}From now on we denote for brevity $\hat U(\xi) \equiv U^{(A)}(\xi)$
and the corresponding $\tau$-function will be
$\hat\tau(\xi,\bar\xi|g)$.\footnote{KP/Toda system itself arises in the limit
$N\longrightarrow \infty$, and $n$ plays the role of the "zero-time" $T_0$ of
the Toda system. In fact, most of equations and their interesting properties
are essentially independent of $N$.} This $\tau$-function is a linear in each
time-variable $\xi_i$, hence, it satisfies simpler determinant formulas and
simpler hierarchy of equations.  Indeed, (\ref{clatau1}) turns into:
\be
\hat\tau_1(\xi,\bar\xi|g) \equiv <0_{F_1}|\hat U(\xi) g \hat{\bar U}(\bar
\xi)|0_{F_1}> = \sum_{k,\bar k \geq 0} s_k\bar s_{\bar k} <k|g|\bar k>
\label{clatau1s}
\ee
where $s_k = \xi_1\xi_2\ldots \xi_{k}$, $s_0=1$, while
(\ref{clataumm1}) is substituted by:
\be
\hat\tau_1^{m\bar m}(\xi,\bar\xi|g)
\equiv <m_{F_1}|\hat U(\xi) g \hat{\bar U}(\bar \xi)|\bar m_{F_1}> =
\frac{1}{s_m\bar s_{\bar m}} \sum_{{k \geq m}\atop{\bar k \geq \bar m}} s_k\bar
s_{\bar k} <k|g|\bar k>   = \\ =\frac{1}{s_{m}\bar s_{\bar m}}\sum_{{k \geq
m}\atop{\bar k \geq \bar m}} \frac{\partial}{\partial\log s_k}
\frac{\partial}{\partial\log \bar s_{\bar k}}\tau_1(\xi,\bar\xi|g)=
\frac{1}{s_{m-1}\bar s_{\bar
m-1}}\frac{\partial}{\partial \xi_m} \frac{\partial}{\partial \bar
\xi_{\bar m}}\tau_1(\xi,\bar\xi|g).
\label{clataumm1s}
\ee

Thus,
\be
\hat\tau_{n+1} = \det_{0\leq m,\bar m \leq n}\hat\tau_1^{m\bar m}
= \left(\prod_{m=1}^n s_m\bar s_{\bar m}\right)^{-1}
\det_{(m,\bar m)} \left(
\sum_{{k \geq m}\atop{\bar k \geq \bar m}} s_k\bar s_{\bar k}
<k|g|\bar k>\right) = \\
=  \frac{1}{s_n\bar s_n}
\sum_{k,\bar k \geq n} s_k\bar s_{\bar k}
\det_{0 \leq m,\bar m \leq n-1}\left(
\begin{array}{cc}
g_{m\bar m}  g_{m\bar k} \\
g_{k\bar m}  g_{m\bar m}
\end{array} \right)
\equiv \frac{1}{s_n\bar s_n}
\sum_{k,\bar k \geq n} s_k\bar s_{\bar k}
{\cal D}_{k\bar k}^{(n)}.
\label{detA}
\ee
One can compare determinant representations (\ref{detB}) and (\ref{detA})
to get the connection between different coordinates $t$ and $\xi$. One can see
that this is of the type $s_k\sim$ some functions of $P_j(t)$.
\footnote{As an example of what it might be like, one can consider the
simplest case of the first fundamental representation. Then, one puts
$\tau_1(t|g) =\hat\tau_1(\xi|g)$
and sees that
$
s_k = P_k(t),\ \ \ {\partial\over\partial t_k}=\sum_i s_{i-k}{\partial\over
\partial s_i}.
$
However, identification of $\tau_n(t)$ and $\hat\tau_n(\xi)$
with $n\neq 1$ will lead to different relations between
$\xi$ and $t$.}

Equations for the $\tau$-function in parameterization A can be easily
derived. Indeed, it is straightforward to find
the Baker-Akhiezer functions (\ref{BA}) and substitute this into
equations (\ref{BIA})
\be
\Psi^{+,n+k+1}_n(\xi)={s_{n+k}\over
s_{n}}\left(\tau_n(\xi)-\xi_n{\partial\tau_n(\xi)\over\partial\xi_n}\right), \\
\Psi^{+,n+1}_n(\xi)=
\left(\tau_n(\xi)-\xi_n{\partial\tau_n(\xi)\over\partial\xi_n}\right),\ \
\Psi^{+,n}_n(\xi)=-{\partial\tau_n(\xi)\over\partial \xi_n};\\
\Psi^{-,k}_n(\xi)-\xi_{n-1}{\partial\Psi^{-,k}_n(\xi)\over\xi_{n-1}}
={s_{n-1}\over s_{k-1}}{\partial\tau_n(\xi)\over\partial
\log\xi_{k}}+{s_{n-1}\over s_{k-2}}{\partial\tau_n(\xi)\over\partial\xi_{k-1}}
\ \hbox{for}\ k>n,\\
\Psi^{-,n}_n(\xi)-\xi_{n-1}{\partial\Psi^{-,n}_n(\xi)\over\xi_{n-1}}
=\tau_n(\xi)+ {\partial\tau_n(\xi)\over \log\xi_{n}},\\
\Psi^{-,n-1}_n(\xi)=\xi_{n-1}\tau_n(\xi),\ \ \Psi^{-,k}_n(\xi)=0\ \hbox{for}\
k<n-1.\label{exBAA}
\ee
As for the values of $\Psi^{+,k}_n(\xi)$ for $k<n$,
they are constants which can be hardly expressed as an action of a differential
operator on $\tau_n(\xi)$. This means that the relation (\ref{BIA}), where
manifest expressions for $\bar\Psi$, analogous to (\ref{exBAA}), can be easily
written down, does not lead to differential equations when $k$ and $l$ are
arbitrarily chosen.  If, however, one chooses $k\le l-1$, because of multiple
cancelations due to (\ref{exBAA}), (\ref{BIA}) is almost a differential
equation.
It can be easily transformed to a differential equation by putting
$\xi_{n-1}=0$ (see (\ref{exBAA})). One can easily check that the
number of independent equations obtained in this way is sufficient to determine
$\tau$-function in full.

\section{Quantum ($q\neq 1$) case}
\subsection{$q$-Determinant-like representation}
In this section we demonstrate how the technique developed in the previous
sections is deformed to the quantum case and, in particular, obtain
$q$-determinant-like representation analogous to (\ref{det}). We also
demonstrate that in parametrization A
relation (\ref{clataumm1s}) expressing $\tau_1^{m\bar m}$ through
$\tau_1$ derivations is still correct for $q\ne 1$, with all the
derivatives replaced by difference operators.

In this subsection we present the statements valid for {\it any} $U(\xi)$ of
the form (\ref{1.11}), without reference to particular parameterization
A.\footnote{Actually, we only require that $U(\xi)$ is an element from $NG_q$
and is expressed only through the generators associated with
{\it simple positive} roots: $U(\xi)=U\{\xi_s|T_i\}$. Formula (\ref{1.11}) is
a possible but not the unique realization of these requirements.}

As a result of the absence of diagonal factor $g_D$
co-product (\ref{coprodcl}) is replaced by the following comultiplication rule
\be
\Delta^{n-1}(U\{T_i\})
= \prod_{m=1}^n U^{(m)}
\ee
where
\be
U^{(m)} = U\left\{\ I\otimes\ldots
\ldots \otimes I\otimes T_i \otimes q^{-2H_i}\otimes \ldots\otimes
q^{-2H_i}\right\}
\ee
($T_i$ appears at the $m$-th place in the tensor
product).  Similarly
\be
\bar U^{(m)} = \bar U
\left\{\ q^{2H_i}\otimes\ldots \ldots \otimes q^{2H_i}\otimes T_{-i}
\otimes I\otimes \ldots\otimes I\right\}.
\ee

Let
$$
H_i|\bar j_{F_1}> = h_{i,\bar j}|\bar
j_{F_1}>, \ \ \ <j_{F_1}|H_i = h_{i,j}<j_{F_1}|
$$
(in fact for $SL(N)$ $\
2h_{i,i-1} = +1,\ 2h_{i,i} = -1$, all the rest are vanishing).  Then
\be
\tau_n^{j_1\ldots j_n\bar j_1\ldots\bar j_n}(\xi_s,\bar \xi_s|g)
\equiv \\ \equiv
\left(\otimes_{m=1}^n<j_m|\right)
\Delta^{n-1}(U)\ g^{\otimes n}\ \Delta^{n-1}(\bar U)
\left(\otimes_{m=1}^n|\bar j_m>\right)=
\\= \prod_{m=1}^n \ <j_m|
U\left\{T_i
q^{-2\sum_{l = m+1}^nh_{i,j_l}}\right\}\ g\
\bar U\left\{T_{-i}
q^{2\sum_{l=1}^{m-1}h_{i,\bar j_l}}\right\}\ |\bar j_m>=\\
=\prod_{m=1}^n \tau_1^{j_m\bar j_m}\left(\xi_s
q^{-2\sum_{l = m+1}^nh_{i(s),j_l}},\bar \xi_s
q^{2\sum_{l=1}^{m-1}h_{i(s),\bar j_l}}\right).
\label{4.42}
\ee

In order to get a $q$-determinant-like counterpart of (\ref{det}),
one should replace  antisymmetrization
by $q$-antisymmetrization in eqs. (\ref{as})-(\ref{ass}), since, in quantum
case, fundamental representations are described by $q$-antisymmetrized vectors
(see s.5.2 of \cite{GKLMM} for more details).
We define $q$-antisymmetrization as a sum over all permutations,
\be
\left( [1,\ldots,k]_q\right) = \sum_P (-q)^{{\rm deg}\ P}
\left(P(1),\ldots,P(k)\right),
\ee
where
\be
{\rm deg}\ P = \#\ {\rm of\ inversions\ in}\ P.
\ee
Then, $q$-antisymmetrizing (\ref{4.42}) with $j_k=k-1,\ \bar j_{\bar
k}=\bar k-1$, one finally gets

\be
\tau_n(\xi,\bar \xi|g) = \sum_{P,P'} (-q)^{{\rm deg}\ P + {\rm deg}\
P'}\times\\
\times\prod_{m=0}^{n-1}
\tau_1^{P(m)P'(\bar m)}\left(\xi_s
q^{-2\sum_{l = m+1}^{n-1}h_{i(s),P(l)}},\bar \xi_s
q^{2\sum_{\bar l=0}^{m-1}h_{i(s),P'(\bar l)}}\right).
\label{detq}
\ee
This would be just a $q$-determinant\footnote{Let us note that the relevant
$q$-determinant is defined as \cite{GKLMM}
\be
{\rm det}_q A \sim  A^{[1}_{[1}\ldots A^{n]_q}_{n]_q}
= \sum_{P,P'} (-q)^{{\rm deg}\ P + {\rm deg}\ P'}
\prod_{a} A^{P(a)}_{P'(a)}.
\label{qdet}
\ee
This is not
necessarily the same as
$A^1_{[1}\ldots
A^n_{n]_q}$. It is the same only provided by special commutation relations of
the matrix elements $A_i^j$.}:
be there no the $q$-factors, which twist the time variables in (\ref{detq}).

To make this expression more transparent, let us consider the simplest
example of the second fundamental representation. Denote through $\{u\}$ and
$\{v\}$ the subsets of $\{s\}$ such that
$i(s)=1$ and $i(s)=2$ respectively. Then
\be\label{tau2}
\tau_2  = \tau_1^{00}(\{q\xi_u\},\{q^{-1}\xi_v\},\xi_s;\
\{\bar \xi_u\},\{\bar \xi_v\},\bar \xi_s) \tau_1^{11}(\{\xi_u\},
\{\xi_v\},\xi_s;\
  \{q\bar \xi_u\},\{\bar \xi_v\},\bar \xi_s) - \\
- q\tau_1^{01}(\{q\xi_u\},\{q^{-1}\xi_v\},\xi_s;\
  \{\bar \xi_u\},\{\bar \xi_v\},\bar \xi_s)
\tau_1^{10}(\{\xi_u\},\{\xi_v\},\xi_s;\
  \{q^{-1}\bar \xi_u\},\{q\bar \xi_v\},\bar \xi_s) - \\
- q\tau_1^{10}(\{q^{-1}\xi_u\},\{\xi_v\},\xi_s;\
  \{\bar \xi_u\},\{\bar \xi_v\},\bar \xi_s)
\tau_1^{01}(\{\xi_u\},\{\xi_v\},\xi_s;\
  \{q\bar \xi_u\},\{\bar \xi_v\},\bar \xi_s) + \\
+ q^2\tau_1^{11}(\{q^{-1}\xi_u\},\{\xi_v\},\xi_s;\
  \{\bar \xi_u\},\{\bar \xi_v\},\bar \xi_s)
\tau_1^{00}(\{\xi_u\},\{\xi_v\},\xi_s;\
  \{q^{-1}\bar \xi_u\},\{q\bar \xi_v\},\bar \xi_s).
\ee
$\xi_s$ here denotes all the time variables with $i(s)>2$.
Let us note that $q$-factors in all these expressions can be reproduced by
action of the operators
$$
M_j^{\pm}:\ \ M_j^{\pm}\xi_{s} =
q^{\pm\delta_{j,i(s)}}\xi_{s},
$$
$$
\bar M_j^{\pm}:\ \ \bar M_j^{\pm}\bar \xi_{s} =
q^{\pm\delta_{j,i(s)}}\bar \xi_{s}.
$$

Now we briefly discuss the set of equations satisfied by quantum
$\tau$-function. We follow the same line as in classical case and introduce
intertwining operators. In quantum case, one should distinguish between the
right and left intertwiners: $\Phi^{\pm,R} :\ \ F_n\otimes F_1$~{\ams
\symbol{29}}~$F_{n+1}$ and $\Phi^{\pm,L} :\ \ F_1\otimes F_n$~{\ams
\symbol{29}}~$F_{n+1}$. These operators $\Phi^{\pm,R,L}$
can be expressed through the classical intertwining operators (fermions):

\be
\Phi^{\pm,R}_i=q^{-\sum_{j=1}^{i-1}\psi^+_j\psi^-_j}\psi^{\pm}_i,\ \ \
\Phi^{\pm,L}_i=q^{\sum_{j=1}^{i-1}\psi^+_j\psi^-_j}\psi^{\pm}_i.
\label{phipsi}
\ee
In analogy with the classical case, one can consider the operator
$\Gamma=\sum_i\Phi^{+,L}_i \otimes\Phi^{-,R}_i$ that
commutes with $g\otimes g$. Then, introducing vertex operators, or
Baker-Akhiezer functions as averages of quantum intertwining operators
(properly labeled by indices $L$ and $R$), one obtains equations
(\ref{BIglinfty})-(\ref{IntHir}) with the properly defined entries.
Technically,
vertex operators can be calculated with the help of equation
(\ref{phipsi}).  In the next subsection we show how it works in the
concrete case of parameterization A.

\subsection{Parameterization A}

Now we apply formulas of the previous subsection to the case
of parameterization A. In fact, most of expressions from subsection
3.3 remain almost the same in the quantum case. In particular,

\be
\hat\tau_1(\xi,\bar\xi|g) \equiv\ <0_{F_1}|\hat U(\xi) g \hat{\bar U}(\bar
\xi)|0_{F_1}>\ = \sum_{k,\bar k \geq 0} s_k\bar s_{\bar k} <k|g|\bar k>
\label{clatau1sq}
\ee
where again $s_k = \xi_1\xi_2\ldots \xi_{k}$, $s_0=1$, while
$\bar s_k=\bar\xi_k\ldots\bar\xi_2\bar\xi_1$, $\bar s_0=1$ and
\be
\hat\tau_1^{m\bar m}(\xi,\bar\xi|g)=
s_m^{-1}
\left(\sum_{{k \geq m}\atop{\bar k \geq \bar m}}
s_k\bar s_{\bar k} <k|g|\bar k>\right) \bar s_{\bar m}^{-1}=\\=
s_{m-1}^{-1}\left( D_{\xi_m} \bar D_{\bar \xi_{\bar m}}\tau_1(\xi,\bar\xi|g)
\right)\bar s_{\bar m-1}^{-1}.
\label{clataumm1sq}
\ee
Here\footnote{There is an ambiguity in the choice of these operators as
$\tau$-function is a linear function of times, and therefore, any linear
operator which makes unity from $\xi$ is suitable. We fix them to
act naturally on the corresponding $q$-exponentials in
(\ref{genpar}).} $ D_{\xi_i} f(\xi)\equiv
\frac{1}{\xi_i}\frac{M_i^{+2}-1}{q^2-1}f(\xi)$,
$ \bar D_{\bar\xi_i}f(\bar\xi) \equiv
\left[\frac{M_i^{-2}-1}{q^{-2}-1}f(\bar\xi)\right]\frac{1}{\bar\xi_i}$.
Then, one can express $\tau_n$ through
$\tau_1$ manifestly using formulas (\ref{detq}) and (\ref{clataumm1sq}).

Equation (\ref{detq}) remains just the same. In our example
(\ref{tau2}) of the second fundamental representation each set
$\{u\}$ and $\{v\}$ consists of the single element:
  $\{u\}=\{s=1\}$, $\{v\}=\{s=2\}$. Then
\be
  \tau_2  = \tau_1^{00}(q\xi_1,q^{-1}\xi_2,\xi_i;\ \bar
\xi_1,\bar \xi_2,\bar \xi_i) \tau_1^{11}(\xi_1,\xi_2,\xi_i;\ q\bar \xi_1,\bar
  \xi_2,\bar \xi_i) - \\ - q\tau_1^{01}(q\xi_1,q^{-1}\xi_2,\xi_i;\ \bar
  \xi_1,\bar \xi_2,\bar \xi_i) \tau_1^{10}(\xi_1,\xi_2,\xi_i;\ q^{-1}\bar
  \xi_1,q\bar \xi_2,\bar \xi_i) - \\ - q\tau_1^{10}(q^{-1}\xi_1,\xi_2,\xi_i;\
  \bar \xi_1,\bar \xi_2,\bar \xi_i)
\tau_1^{01}(\xi_1,\xi_2,\xi_i;\
  q\bar \xi_1,\bar \xi_2,\bar \xi_i) + \\
+ q^2\tau_1^{11}(q^{-1}\xi_1,\xi_2,\xi_i;\
  \bar \xi_1,\bar \xi_2,\bar \xi_i)
\tau_1^{00}(\xi_1,\xi_2,\xi_i;\
  q^{-1}\bar \xi_1,q\bar \xi_2,\bar \xi_i).
\ee
This expression can be written in a more compact form with
the help of operators
\be
{\cal D}_1^L \equiv M_1^-D_1\otimes I, \ \ \ {\cal D}_1^R
\equiv M_1^+M_2^-\otimes D_1, \\ \bar{\cal D}_1^L \equiv
{\bar D_1} \otimes \bar M_1^-\bar M_2^+, \ \ \ \bar{\cal
D}_1^R \equiv I\otimes \bar M_1^+{\bar D_1}.
\ee
Indeed,
\be
\tau_2 = \left({\cal D}_1^R\bar{\cal D}_1^R -
q{\cal D}_1^L\bar{\cal D}_1^R - q{\cal D}_1^R\bar{\cal D}_1^L
+ q^2{\cal D}_1^L\bar{\cal D}_1^L\right)\tau_1\otimes
\tau_1 =\\=\left( {\cal D}_1^R-q{\cal D}_1^L\right)\cdot\left(
\bar{\cal D}_1^R-q\bar{\cal D}_1^L\right)
\tau_1\otimes \tau_1.
\ee
These operators satisfy the following commutation
relations (like the algebra of $\theta$ and $\chi$ in (\ref{comrel})):
\be
{\cal D}_1^{L}    {\cal D}_1^{R}
= q{\cal D}_1^{R}{\cal D}_1^{L}, \\
\bar {\cal D}_1^{L} \bar   {\cal D}_1^{R}
= q\bar{\cal D}_1^{R}\bar{\cal D}_1^{L}.
\ee
These formulas can be rewritten in a more "invariant" form
in terms of operators
\be
\hbox{{\hr D}}^L_i\equiv D_i\otimes I,\ \ \
\hbox{{\hr D}}^R_i\equiv \prod_j M_j^{-\vec\alpha_i\vec\alpha_j}\otimes D_i,
\\
\bar {\hbox{{\hr D}}}^L_i\equiv \bar D_i\otimes\prod_j \bar
M_j^{-\vec\alpha_i\vec\alpha_j},\ \ \
\bar {\hbox{{\hr D}}}^R_i\equiv I\otimes \bar D_i,
\ee
which commute as
\be
\hbox{{\hr D}}_i^{L}    \hbox{{\hr D}}_j^{R}
= q^{\vec\alpha_i\vec\alpha_j}
\hbox{{\hr D}}_j^{R}\hbox{{\hr D}}_i^{L}, \\
\bar {\hbox{{\hr D}}}_i^{L} \bar {\hbox{{\hr D}}}_j^{R}
= q^{\vec\alpha_i\vec\alpha_j}\bar{\hbox{{\hr D}}}_j^{R}\bar{\hbox{{\hr
D}}}_i^{L}.
\ee
Then,
\be
\tau_2 = M^-_1\otimes\bar M_1^+\left( \hbox{{\hr D}}_1^R-q
\hbox{{\hr D}}_1^L\right)\cdot\left(
\bar{\hbox{{\hr D}}}_1^R-q\bar{\hbox{{\hr D}}}_1^L\right)
\tau_1\otimes \tau_1.
\ee

Baker-Akhiezer functions for the $\tau$-function in
parametrization A are given by the following expressions
\be
\Psi^{+,n+k+1}_n(\xi)=q^{n+1}s_n^{-1}s_{n+k}M^+_{n+1}\ldots M^+_{n+k-1}
\left(\tau_n(\xi)-\xi_nD_n\tau_n(\xi)\right), \\
\Psi^{+,n+1}_n(\xi)=q^{n+1}
\left(\tau_n(\xi)-\xi_nD_n\tau_n(\xi)\right),\ \
\Psi^{+,n}_n(\xi)=-q^nD_n\tau_n(\xi);\\
\Psi^{-,k}_n(\xi)-\xi_{n-1}D_{n-1}\Psi^{-,k}_n(\xi)=\\=q^{n-2}
s_{k-2}^{-1}s_{n-1}D_{i-1}\tau_n(\xi)+q^{n-2}
s_{k-1}^{-1}s_{n-1}\xi_iD_i\tau_n(\xi)
\ \hbox{for}\ k>n,\\
\Psi^{-,n}_n(\xi)-\xi_{n-1}D_{n-1}\Psi^{-,n}_n(\xi)=q^{n-2}
\tau_n(\xi)+ \xi_nD_n\tau_n(\xi),\\
\Psi^{-,n-1}_n(\xi)=q^{n-2}\xi_{n-1}\tau_n(\xi),\ \ \Psi^{-,k}_n(\xi)=0\
\hbox{for}\ k<n-1.
\ee
Substituting these expressions to
formula (\ref{BIA}), one obtains the set of equations which is a quantum
counterpart of KP/Toda hierarchy in parameterization A.

\section{Conclusion}
In the paper we described the general way to construct quantum deformations of
the determinant representations of KP/Toda $\tau$-functions. For this, we
  did not need any concrete form of the time evolution but only the suggestion
  that the evolution is described by (quasi-)group elements.

We observed that,
for $q\ne 1$, the determinant representations turn into
$q$-determinant-like ones. Moreover, we do not get just
$q$-determinants only because the evolution operator in quantum case is not a
group element. This happens,
because no nilpotent subgroup $NG_q$ exists in the quantum group. To avoid this
  problem, one could begin from a slightly different parametrization of
classical $\tau$-function, such that the evolution operator lies in Borel (not
just nilpotent) subgroup $BG_q$. In the classical limit, additional Cartan
generators can be removed by redefinition of the element $g$ labeling the point
of the Grassmannian,but, in the quantum case, the Cartan part of the
evolution would essentially change the result: evolution operator is now a
group element (for $BG_q$) and, thus, additional twistings of times disappear
from formula (\ref{detq}), and so defined quantum $\tau$-functions are just the
$q$-determinants. In order to fulfil this program, one still needs to find an
appropriate parametrization of (a set of) the group elements of $BG_q$ by
exactly $r_G$ "time-variables". In this paper we followed another way:
evolution has been taken to lie in $NG_q$, and we explicitly described the
corrections to naive $q$-deformed formulas, which originated from the fact that
$NG_q$ is not a subgroup.

Another problem, which remains beyond the scope of this paper, is deformation
of
the standard KP/Toda evolutions B and C. We mentioned in section 2 that this
problem is equivalent to restricting
possible representations of algebra (\ref{comrel}) of functions ($\theta,\
\chi,\ \phi$) in such a way that the number of independent
variables is appropriately reduced. Unfortunately, we are not aware of
explicit solutions to these constraints neither in the B, nor in
C case.

The third problem, which deserves to be mentioned here, concerns
interpretation of integrable hierarchies in terms of the Grassmannian. In fact,
in this paper we derived hierarchies of integrable equations making use of
intertwining operators (fermions). These latter ones give a natural
definition of $q$-deformed fermions, each fermion being doubled due to
the difference between right and left intertwiners in quantum group.
We observed that, in order to derive the equations, one needed bilinear
combinations mixing the right and the left
intertwiners\footnote{Generators of quantum enveloping algebra also can be
described by the same "mixed" combinations.}.  On the other hand, in order to
properly describe the $q$-Grassmannian, one might need combinations of
{\it all} intertwiners including those with only left, or only right
ones. In this case, in order to obtain an adequate description of quantum
hierarchies, one needs a sort of a doubled $q$-Grassmannian.

The details of description of quantum hierarchies in terms of quantum
intertwining operators ($q$-fermions) and related problems will be presented in
\cite{Det}.

\section*{Acknowledgements}
We are grateful to A.Gerasimov for useful discussions.
The work of S.Kh. and A.Mir. is partially supported by grant 93-02-03379 of the
Russian Foundation of Fundamental Research and by INTAS grant 93-1038.

\end{document}